\documentclass{amsart}
\usepackage{amsmath, amssymb, latexsym, amsthm, mathrsfs}

\newcommand{\Pa}[9]{\bibitem{#1} {#2}, \emph{#3}, {#4} \textbf{#5} ({#6}), {#7}--{#8}.}

\newcommand{\ed}{

\end{document}
}

\newcommand{\Arh}{Arhangel'ski\u{\i}}

\newcommand{\Setting}[7]{\xymatrix@R=4pt@C=7pt{#1\ar@{-}[r]&#2\ar@{-}[r]&#3\\&#4\ar@{-}[u]\\
#5\ar@{-}[uu]\ar@{-}[r] & #6\ar@{-}[u]\ar@{-}[r] & #7\ar@{-}[uu]}}

\newcommand{\N}{\mathbb{N}}

\newcommand{\compactN}{\cl{\mathbb{N}}}

\newcommand{\arx}[3]{\texttt{http://arxiv.org/#1/#3}}
\newcommand{\bq}{\begin{quote}}
\newcommand{\eq}{\end{quote}}
\newcommand{\cl}[1]{\overline{#1}}
\newcommand{\CH}{the Continuum Hypothesis}
\newcommand{\MA}{Martin's Axiom}

\newcommand{\NN}{{\N^{\N}}}

\newcommand{\Inc}{{\compactN^{\uparrow\N}}}

\newcommand{\NcompactN}{{\compactN^\N}}
\newcommand{\seq}[1]{\{#1\}_{n\in\N}}
\newcommand{\sseq}[1]{\{#1 : n\in\N\}}

\newcommand{\op}{\operatorname}

\newcommand{\scrA}{\mathscr{A}}
\newcommand{\scrB}{\mathscr{B}}

\newcommand{\CG}{C_\Gamma}

\newcommand{\cF}{\mathcal{F}}

\newcommand{\cM}{\mathcal{M}}

\newcommand{\cO}{\mathcal{O}}

\newcommand{\R}{\mathbb{R}}
\newcommand{\cU}{\mathcal{U}}
\newcommand{\Union}{\bigcup}

\long\def\forget#1\forgotten{}

\newcommand{\fb}{\mathfrak{b}}
\newcommand{\fc}{\mathfrak{c}}
\newcommand{\fd}{\mathfrak{d}}

\newcommand{\oo}{\infty}

\newcommand{\nin}{\notin}
\newcommand{\cat}{\hat{\ }}
\newcommand{\sbst}{\subseteq}
\newcommand{\spst}{\supseteq}
\newcommand{\sm}{\setminus}

\newcommand{\rest}{\restriction}

\newcommand{\cov}{\op{cov}}

\newtheorem{thm}{Theorem}[section]
\newtheorem{prop}[thm]{Proposition}

\newtheorem{prob}[thm]{Problem}
\newtheorem{lem}[thm]{Lemma}
\newtheorem{cor}[thm]{Corollary}

\theoremstyle{definition}

\theoremstyle{remark}

\newcommand{\be}{\begin{enumerate}}
\newcommand{\ee}{\end{enumerate}}
\newcommand{\bi}{\begin{itemize}}
\newcommand{\itm}{\item}
\newcommand{\ei}{\end{itemize}}
\forget
\forgotten

\newcommand{\sone}{\mathsf{S}_1}
\newcommand{\sfin}{\mathsf{S}_\mathrm{fin}}
\newcommand{\ufin}{\mathsf{U}_\mathrm{fin}}

\title{Hurewicz sets of reals without perfect subsets}

\author{Du\v{s}an Repov\v{s}}
\address{Du\v{s}an Repov\v{s}, Institute of  Mathematics, Physics and
Mechanics and Faculty of Education, University of Ljubljana,
P.O.B. 2964, Ljubljana, Slovenija 1001.}
\email{dusan.repovs@guest.arnes.si}

\author{Boaz Tsaban}
\address{Boaz Tsaban,
Department of Mathematics, Bar-Ilan University, Ramat-Gan 52900, Israel;
and
Department of Mathematics, Weizmann Institute of Science, Rehovot 76100, Israel.
}
\email{tsaban@math.biu.ac.il}

\author{Lyubomyr Zdomskyy}
\address{Lyubomyr Zdomskyy, Kurt G\"odel Research Center for Mathematical Logic, W\"ahringer Str.\ 25, A-1090 Vienna, Austria.}
\email{lzdomsky@gmail.com}

\thanks{The first and the third authors were supported
by the Slovenian Research Agency grants P1-0292-0101-04 and BI-UA/04-06-007.
The second author was partially supported by the Koshland Center for Basic Research.}


\begin{document}

\begin{abstract}
We show that even for subsets $X$ of the real line which do not contain
perfect sets, the Hurewicz property does not imply the property $\sone(\Gamma,\Gamma)$,
asserting that for each countable family of open $\gamma$-covers of $X$, there is a choice
function whose image is a $\gamma$-cover of $X$.
This settles a problem of Just, Miller, Scheepers, and Szeptycki. Our main result
also answers a question of Bartoszy\'nski and the second author,
and implies that for $C_p(X)$, the conjunction of Sakai's strong countable fan tightness
and the Reznichenko property does not imply \Arh{}'s property $\alpha_2$.
\end{abstract}

\maketitle

\section{Introduction}

By a \emph{set of reals} we mean a separable, zero-dimensional, and metrizable space
(such spaces are homeomorphic to subsets of the real line $\R$).
Fix a set of reals $X$. Let $\cO$ denote the collection of all open covers of $X$.
An open cover $\cU$ of $X$ is a \emph{$\gamma$-cover} of $X$ if it is infinite and for each $x\in X$,
$x$ is a member of all but finitely many members of $\cU$.
Let $\Gamma$ denote the collection of all open $\gamma$-covers of $X$.
Motivated by Menger's work, Hurewicz \cite{Hure27} introduced the \emph{Hurewicz property}
$\ufin(\cO,\Gamma)$:
\bq
For each sequence $\seq{\cU_n}$ of members of $\cO$
which do not contain a finite subcover,
there exist finite sets $\cF_n\sbst\cU_n$, $n\in\N$,
such that $\sseq{\cup\cF_n}\in\Gamma$.
\eq
Every $\sigma$-compact space satisfies $\ufin(\cO,\Gamma)$, but the converse
fails \cite{coc2, ideals}.

Let $\scrA$ and $\scrB$ be any two families.
Motivated by works of Rothberger, Scheepers introduced the following
prototype of properties \cite{coc1}.
\begin{description}
\item[$\sone(\scrA,\scrB)$]
For each sequence $\seq{\cU_n}$ of members of $\scrA$,
there exist members $U_n\in\cU_n$, $n\in\N$, such that $\sseq{U_n}\in\scrB$.
\end{description}
It is easy to see that $\ufin(\cO,\Gamma)=\ufin(\Gamma,\Gamma)$, and
therefore $\sone(\Gamma,\Gamma)$ implies $\ufin(\cO,\Gamma)$ \cite{coc1}.
However, a set of reals satisfying $\sone(\Gamma,\Gamma)$
cannot contain perfect subsets \cite{coc2}. It follows that, for example, $\R$ satisfies
$\ufin(\cO,\Gamma)$ but not $\sone(\Gamma,\Gamma)$.
In the fundamental paper \cite{coc2}, we are asked whether there are \emph{nontrivial}
examples showing that $\ufin(\cO,\Gamma)$ does not imply $\sone(\Gamma,\Gamma)$.

\begin{prob}[Just, Miller, Scheepers, Szeptycki \cite{coc2}]\label{prob}
Let $X$ be a set of reals which does not contain a perfect set,
but which does have the Hurewicz property. Does $X$ then satisfy $\sone(\Gamma,\Gamma)$?
\end{prob}

We give a negative answer that also yields a new result concerning function spaces.

\section{The main theorem}

We prove a stronger assertion than what is needed to settle Problem \ref{prob};
this will be useful for the next section. Let
$\CG$ denote the collection of all \emph{clopen} $\gamma$-covers of $X$.
Clearly, $\sone(\Gamma,\Gamma)$ implies $\sone(\CG,\CG)$.\footnote{It is an open problem whether
the converse implication holds \cite{BH07, Sakai07}.}
The hypothesis in the following theorem is a consequence of \CH{}.
See \cite{BlassHBK} for a survey of the involved cardinals.

\begin{thm} \label{maintheorem}
Assume that $\fb=\fc$.
There exists a set of reals $X$ such that:
\be
\itm $X$ does not contain a perfect set;
\itm All finite powers of $X$ have the Hurewicz property $\ufin(\cO,\Gamma)$; and
\itm No set of reals containing $X$ satisfies $\sone(\CG,\CG)$.
\ee
\end{thm}

Theorem \ref{maintheorem} is proved in three steps.
The first step is analogous to Theorem 4.2 of \cite{BRR91},
and will be used to show that the constructed set is not contained in a set of reals
satisfying $\sone(\CG,\CG)$.
Say that a convergent sequence $\seq{x_n}$ is \emph{nontrivial}
if $\lim_nx_n\nin\sseq{x_n}$.

\begin{lem} \label{lem1}
Let $X$ be a subspace of a zero-dimensional metrizable space $Y$ satisfying $\sone(\CG,\CG)$,
and $\seq{x^m_n}$, $m\in\N$, be nontrivial convergent sequences in $X$.
Then: There are a countable closed cover $\{F_k:k\in\N\}$ of $X$ and an infinite $A\sbst\N$,
such that $F_k\cap\{x^m_n:n\in A\}$ is finite for all $k,m$.
\end{lem}
\begin{proof}
Let $d$ be a metric on $Y$ which generates its topology.
For each $m$, do the following. Let $x_m=\lim_n x^m_n$, and for each $n$ take a clopen neighborhood
$C^m_n$ of $x^m_n$ in $Y$, whose diameter is smaller than $d(x^m_n,x_m)/2$.
For each $m,n$, set
$$U^m_n = Y\sm (C^0_n\cup C^1_n\cup \dots \cup C^m_n).$$
For each $m$, $\sseq{U^m_n}$ is a clopen $\gamma$-cover of $Y$.
Apply $\sone(\CG,\CG)$ to get $f\in\NN$ such that $\{U^m_{f(m)}:m\in\N\}$
is a (clopen) $\gamma$-cover of $Y$.
As $U^m_{f(m)}\sbst Y\sm C^0_{f(m)}$ for each $m$,
we have that the image $A$ of $f$ is infinite.

For each $k$, let $F_k = \bigcap_{i\ge k}U^i_{f(i)}$.
$\{F_k : k\in\N\}$ is a closed ($\gamma$-)cover of $Y$.
Fix $k$ and $m$.
If $n$ is large enough and $n\in A$, then $n=f(i)$ with
$i\ge m,k$.
As $x^m_n=x^m_{f(i)}\in C^m_{f(i)}$
and $i\ge m$, $x^m_n\nin U^i_{f(i)}$.
As $i\ge k$, $U^i_{f(i)}\spst F_k$, and therefore
$x^m_n\nin F_k$.
\end{proof}

To make sure that our constructed set does not contain a perfect set and
that it satisfies the Hurewicz property in all finite
powers, we will use the following.
Let $\compactN=\mathbb{N}\cup\{\infty\}$ be the one point compactification of $\N$, and
$\Inc$ be the collection of all nondecreasing elements $f$ of $\NcompactN$ (endowed with the
Tychonoff product topology) such that $f(n)<f(n+1)$ whenever $f(n)<\infty$.
$\Inc$ is homeomorphic to the Cantor space (see \cite{SFH} for an explicit homeomorphism),
and can therefore be viewed as a set of reals.

Let $S$ be the family of all nondecreasing finite sequences in $\N$.
For $s\in S$, $|s|$ denotes its length.
For each $s\in S$, define $q_{s} \in \Inc$ by
$q_s(n) = s(n)$ if $n <|s|$, and $q_s(n) = \infty$ otherwise.
Let $Q$ be the collection of all these elements $q_s$.
$Q$ is dense in $\Inc$.

For a set $D$ and $f,g\in \N^D$, $f\le^* g$ means: $f(d)\le g(d)$ for all but finitely many $d\in D$.
A \emph{$\fb$-scale} is an unbounded (with respect to $\le^*$) set $\{f_\alpha : \alpha<\fb\}\sbst\NN$
of increasing functions, such that $f_\alpha\le^* f_\beta$ whenever $\alpha<\beta$.

\begin{thm}[Bartoszy\'nski-Tsaban \cite{ideals}]\label{H}
Let $X\sbst\Inc$ be a union of a $\fb$-scale and $Q$.
Then $X$ contains no perfect subset, and all finite powers of $H$ satisfy the Hurewicz property
$\ufin(\cO,\Gamma)$.
\end{thm}

For each $s\in S$, $\seq{q_{s\cat n}}$ (where $\cat$ denotes a concatenation of sequences)
is a nontrivial convergent sequence in $\Inc$,\footnote{Strictly speaking,
$q_{s\cat n}\nin\Inc$ when $n<s(|s|-1)$, but since we are dealing with convergent
sequences, we can ignore the first few elements.}
and
$$\lim_{n\to\infty} q_{s\cat n}=q_s.$$
The following will be used in our construction.

\begin{lem} \label{lem2}
Let $X$ be a closed subspace of $\Inc$. If $X\cap \sseq{q_{s\cat n}}$ is finite
for each $s\in S$, then there exists
$\phi:S\to\N$ such that for all $x\in X$ and all $n\geq 2$,
$x(n)\geq\phi(x\rest n)$ implies $x(n+1)\leq\phi(x\rest (n+1))$.
\end{lem}
\begin{proof}
For each $s\in S$, let $k(s)$ be such that $q_{s\cat k}\in \Inc\sm X$ for all $k\geq k(s)$.
As $X$ is closed in $\Inc$, for each $k\geq k(s)$ there is $m(s,k)$ such that
$$\left \{z\in\Inc : z\rest (|s|+1) = s\cat k,\ z(|s|+1)>m(s,k)\right \}\cap X = \emptyset.$$
(Note that $\{z\in\Inc : z\rest (|s|+1) = s\cat k,\ z(|s|+1)>m\}$, $m\in\N$, is
a neighborhood base at $q_{s\cat k}$.)
Define $\phi:S\to\N$ by
$$\phi(s) = \max\{k(s), m(s\rest (|s|-1)), s(|s|-1)\}$$
when $|s|\ge 2$, and by $\phi(s)=0$ when $|s|<2$.
Let $x\in X$ and $n\geq 2$.
If $x(n)\geq\phi(x\rest n)$, then $x(n)\geq k(x\rest n)$, and hence $x(n+1)\leq m(x\rest n,x(n))\leq\phi(x\rest (n+1))$.
\end{proof}

It remains to prove the following.

\begin{prop} \label{lem3}
Assume that $\fb=\fc$.
There exists a $\fb$-scale $B=\{b_\alpha:\alpha<\fb\}$
such that for each closed cover $\sseq{F_n}$ of $B\cup Q$ and each infinite set $A\sbst\N$,
there are $n$ and $s\in S$ such that $F_n\cap\{q_{s\cat k}:k\in A\}$ is infinite.
\end{prop}
\newcommand{\CaA}{\cl{A}^{\uparrow\N}}
\begin{proof}
Let $\{A_\alpha:\alpha<\fc\}$ be an enumeration of all infinite subsets of $\N$,
such that for each infinite $A\sbst\N$, there are $\fc$ many
$\alpha<\fc$ with $A_\alpha=A$.

As $\fb=\fd=\fc$, there is a (standard) scale in $\N^S$, that is,
a family $\{\phi_\alpha:\alpha<\fc\}\sbst \N^S$ such that:
\be
\itm For each $\phi\in \N^S$, there is $\beta<\fc$
such that $\phi\le^* \phi_\beta$; and
\itm For all $\alpha<\beta<\fc$, $\phi_\alpha\le^* \phi_\beta$.
\ee
For an infinite $A\sbst\N$, let $\cl{A}=A\cup\{\infty\}$, and
$$\CaA=\{x\in\Inc:x(n)\in\cl{A} \mbox{ for all }n\}$$
The order isomorphism between $A\cup\{\oo\}$ and $\N\cup\{\oo\}$
induces an order isomorphism $\Psi_A:\CaA\to\Inc$.

By induction on $\alpha<\fb=\fc$, construct a $\fb$-scale $B=\{b_\alpha:\alpha<\fc\}$
such that for each $\alpha<\fc$, $b_\alpha\in (A_\alpha)^{\uparrow\N}$, and
$$\Psi_{A_\alpha}(b_\alpha)(n)>
\phi_\alpha(\Psi_{A_\alpha}(b_\alpha)\rest n)$$
for all $n\geq 2$.

We claim that $X=B\cup Q$ is as required.
Indeed, let $A$ be an infinite subset of $\N$.
Take an increasing enumeration $\{\beta_\alpha:\alpha<\fc\}$ of $\{\alpha<\fc:A_\alpha=A\}$.
For each $\alpha<\fc$, $b_{\beta_\alpha}\in\CaA$.
Set $c_\alpha=\Psi_A(b_{\beta_\alpha})$,
and $C=\{c_\alpha:\alpha<\fc\}$.
By the construction of the functions $b_\alpha$,
$$c_\alpha(n)>\phi_{\beta_\alpha}(c_\alpha\rest n)\geq\phi_\alpha(c_\alpha\rest n)$$
for all but finitely many $n$.

Let  $\{K_m:m\in\N\}$ be a closed cover of $C\cup Q$.
Then there are
$m$ and $s\in S$ such that $K_m\cap\{q_{s\cat k}:k\in\N\}$ is infinite:
Otherwise,  by Lemma \ref{lem2}, for each $m$ there is
$\psi_m\in \N^S$ such that for all $x\in K_m$ and $n\geq 2$,
$x(n)\geq\psi_m(x\rest n)$ implies $x(n+1)\leq\psi_m(x\rest (n+1))$.
Let $\alpha<\fc$ be such that for each $m$, $\phi_\alpha(s)\geq \psi_m(s)$
for all but finitely many $s\in S$. It is easy to verify that
$c_\alpha\not\in K_m$ for all $m$; a contradiction.

Now consider any closed cover $\{F_m:m\in\N\}$ of $B\cup Q$ and set
$K_m=\Psi_A(F_m\cap \CaA)$. Let $s\in S$ and $m$ be such that
$K_m\cap\{q_{s\cat k}:k\in\N\}$ is infinite. Then for $\tilde{s}\in S$
such that $\tilde{s}(i)$ is the $s(i)$'th element of $A$ for each $i<|s|$, we have
that $F_m\cap\{q_{\tilde{s}\cat k}:k\in A\}$ is infinite.
\end{proof}

This completes the proof of Theorem \ref{maintheorem}.
The following corollary of Theorem \ref{maintheorem}
answers in the negative Problem 15(1) of Bartoszy\'nski
and the second author \cite{ideals}.

\begin{cor}
The union of a $\fb$-scale and $Q$ need not satisfy $\sone(\Gamma,\Gamma)$.\hfill\qed
\end{cor}

\section{Reformulation for spaces of continuous functions}

Let $Y$ be a (not necessarily metrizable) topological space.
For $y\in Y$ and $A\sbst Y$, write $\lim A=y$ if $A$ is countable,
and an (any) enumeration of $A$ converges nontrivially to $y$.
Let $\Gamma_y = \{A\sbst Y  : \lim A=y\}$.
$Y$ has the Arhangel'ski\u{\i} \emph{property $\alpha_2$} \cite{Arhan72}
if $\sone(\Gamma_y,\Gamma_y)$ holds for all $y\in Y$.

Fix a set of reals $X$. $C_p(X)$ is the subspace of the Tychonoff
product $\R^X$ consisting of the continuous functions.
It was recently discovered, independently by Bukovsk\'y-Hale\v{s} \cite{BH07} and
by Sakai \cite{Sakai07}, that $C_p(X)$ has the property $\alpha_2$
if, and only if, $X$ satisfies $\sone(\CG,\CG)$.

Many additional connections of this type are studied in the literature.
For families $\scrA$ and $\scrB$, consider the following prototype \cite{coc1}.
\begin{description}
\item[$\sfin(\scrA,\scrB)$]
For each sequence $\seq{\cU_n}$ of members of $\scrA$,
there exist finite subsets $\cF_n\sbst\cU_n$, $n\in\N$, such that $\Union_n\cF_n\in\scrB$.
\end{description}
For a topological space $Y$ and $y\in Y$, let
$\Omega_y = \{A\sbst Y : y\in\cl{A}\sm A\}$.
$Y$ has the Arhangel'ski\u{\i} \emph{countable fan tightness} \cite{Arhan72} if
$\sfin(\Omega_y,\Omega_y)$ holds for each $y\in Y$.
$Y$ has the \emph{Reznichenko property} if for each $y\in Y$
and each $A\in \Omega_y$, there are pairwise disjoint finite sets $F_n\sbst A$, $n\in\N$,
such that each neighborhood $U$ of $y$ intersects $F_n$ for all but finitely many
$n$.

For sets of reals $X$,
$C_p(X)$ has countable fan tightness and the Reznichenko property
if, and only if, all finite powers of $X$ have the Hurewicz property $\ufin(\cO,\Gamma)$
\cite{coc7}.
Thus, Theorem \ref{maintheorem} can be reformulated as follows.

\begin{thm}
Assume that $\fb=\fc$. There exists a set of reals $X$ without perfect subsets,
such that $C_p(X)$ has countable fan tightness and the Reznichenko property,
but does not have the \Arh{} property $\alpha_2$.\hfill\qed
\end{thm}

A topological space $Y$ has the \emph{Sakai strong countable fan tightness}
if $\sone(\Omega_y,\Omega_y)$ holds for each $y\in Y$.
Sakai proved that for sets of reals, $C_p(X)$ has strong countable fan tightness
if, and only if, all finite powers of $X$ satisfy $\sone(\cO,\cO)$ \cite{Sakai88}.
For sets of reals $X$,
$C_p(X)$ has strong countable fan tightness and the Reznichenko property
if, and only if, all finite powers of $X$ satisfy $\ufin(\cO,\Gamma)$ as well as $\sone(\cO,\cO)$
\cite{KocSch02}.

If $\fb\le\cov(\cM)$ and $X$ is a union of a $\fb$-scale and $Q$, then
all finite powers of $X$ satisfy $\ufin(\cO,\Gamma)$ as well as $\sone(\cO,\cO)$ \cite{ideals}.
As \CH{} (or just \MA{}) implies that $\fb=\cov(\cM)=\fc$, we have the following.

\begin{cor}
Even for $C_p(X)$ where $X$ is a set of reals,
the conjunction of strong countable fan tightness and the Reznichenko property
does not imply the \Arh{} property $\alpha_2$.\hfill\qed
\end{cor}

\section{Concluding remarks and open problems}

Our results are consistency results. What is not settled is whether the answers
to the problems addressed in this paper are undecidable.

\begin{prob}
Is it consistent that all sets of reals which have the Hurewicz property $\ufin(\cO,\Gamma)$
but have no perfect subsets satisfy $\sone(\Gamma,\Gamma)$?
\end{prob}

\begin{prob}
Is it consistent that each union of a $\fb$-scale and $Q$ satisfies:
\be
\itm $\sone(\Gamma,\Gamma)$?
\itm $\sone(\Gamma,\Gamma)$ in all finite powers?
\ee
\end{prob}

\begin{prob}
Is it consistent that for each set of reals $X$,
if $C_p(X)$ has both strong countable fan tightness and the Reznichenko property,
then $C_p(X)$ has the \Arh{} property $\alpha_2$?
\end{prob}

\ed